\documentclass[11pt]{amsart}

\usepackage{graphicx,amssymb,epstopdf,subcaption,mathtools,tikz-cd,pb-diagram,thmtools, thm-restate,pinlabel,pifont,enumitem,stmaryrd}
\usepackage{pinlabel}
\usepackage[hmargin=1.55in, vmargin=1.30in]{geometry}
\usepackage[]{hyperref}
\usepackage{tikz}
\usetikzlibrary{matrix}

\usepackage[all]{xy}
\usepackage[mathscr]{euscript}

\newcommand{\nc}{\newcommand}
\nc{\dmo}{\DeclareMathOperator}
\nc{\nt}{\newtheorem}
\nc{\p}[1]{\bigskip\noindent\emph{#1.}}

\nt{theorem}{Theorem}[section]
\nt{corollary}[theorem]{Corollary}
\nt{lemma}[theorem]{Lemma}
\nt{defn}[theorem]{Definition}
\nt{criterion}[theorem]{Criterion}

\nt{question}[theorem]{Question}
\nt{problem}[question]{Problem}
\nt{conjecture}[question]{Conjecture}
\nt{proposition}[theorem]{Proposition}\nt*{claim}{Claim}

\nc{\M}{\mathcal{M}}
\nc{\X}{\mathcal{X}}
\nc{\C}{\mathcal{C}}
\nc{\T}{\mathcal{T}}
\nc{\W}{\mathcal{W}}
\nc{\J}{\mathcal{K}}
\nc{\ch}{e}

\nc{\cut}{\!\ssearrow\!}

\dmo{\Mod}{Mod}
\dmo{\SMod}{SMod}
\dmo{\I}{\mathcal{I}}
\dmo{\SL}{SL}
\dmo{\PSp}{PSp}
\dmo{\PSL}{PSL}
\dmo{\Homeo}{Homeo}
\dmo{\im}{\mathrm{im}}
\dmo{\Aut}{Aut}
\dmo{\Ch}{Ch}
\dmo{\Sp}{\mathrm{Sp}({2g}, \Z)}
\dmo{\lcm}{lcm}

\nc{\Z}{\mathbb Z}
\nc{\N}{\mathcal N}
\nc{\R}{\mathbb R}
\nc{\F}{\mathcal F}

\nc{\ga}{\gamma}
\nc{\de}{\delta}
\nc{\ep}{\epsilon}

\nc{\flm}{\lambda_{2}}

\nc{\normalclosure}[1]{\ensuremath{\left \langle \left \langle #1 \right \rangle \right \rangle}}

\nc{\margin}[1]{\marginpar{\scriptsize #1}}

\title{The normal closure of a bounding pair map in its mapping class group}
\author{Lei Chen}
\author{Weiyan Chen}
\author{Justin Lanier}

\address{Lei Chen \newline Morningside Center of Mathematics, Chinese Academy of Sciences \newline  Academy of Mathematics and Systems Science, Chinese Academy of Sciences\newline  Beijing, 100190, China \\  chenlei@amss.ac.cn }
\address{Weiyan Chen \newline Yau Mathematical Sciences Center, Tsinghua University\newline Beijing, 100084, China \\  chwy@tsinghua.edu.cn}
\address{Justin Lanier \newline Department of Mathematics, Louisiana State University\newline Baton Rouge, LA, 70802, USA \\  justin.lanier@lsu.edu}

\begin{document}

\begin{abstract}
Johnson showed that the normal subgroup of a mapping class group generated by the genus 1 bounding pair maps is equal to the Torelli group. Generalizing Johnson's result, we give two descriptions of the normal subgroup generated by the genus $n$ bounding pair maps using the Chillingworth and the Casson--Morita invariants. 
\end{abstract}

\maketitle

\vspace*{-4ex}



\vspace*{0in}
\vspace{.15in}

\vspace{-.15in}

\section{Introduction}

Birman gave the first generating sets for the Torelli groups of closed surfaces as a consequence of producing a finite presentation for $\Sp$ \cite{BirmanGen}. Powell realized Birman's generators geometrically and showed that the Torelli group of a closed surface is generated by genus $1$ bounding pair maps along with genus $1$ and genus $2$ separating twists \cite{Powell}. Johnson shortly thereafter trimmed down this generating set, showing that for closed and also once-bordered surfaces of genus at least 3, the genus $1$ bounding pair maps suffice \cite{JohnsonLantern}. Given Johnson's result, a natural question presents itself: what subgroup of a Torelli group is generated by the genus $n$ bounding pair maps? Or equivalently, since all bounding pair maps of the same genus are conjugate in their mapping class group: what is the normal closure of a genus $n$ bounding pair map in its mapping class group?

\begin{figure}[ht]
    \centering
        \vspace{-7pt}
    \includegraphics[width=0.45\linewidth]{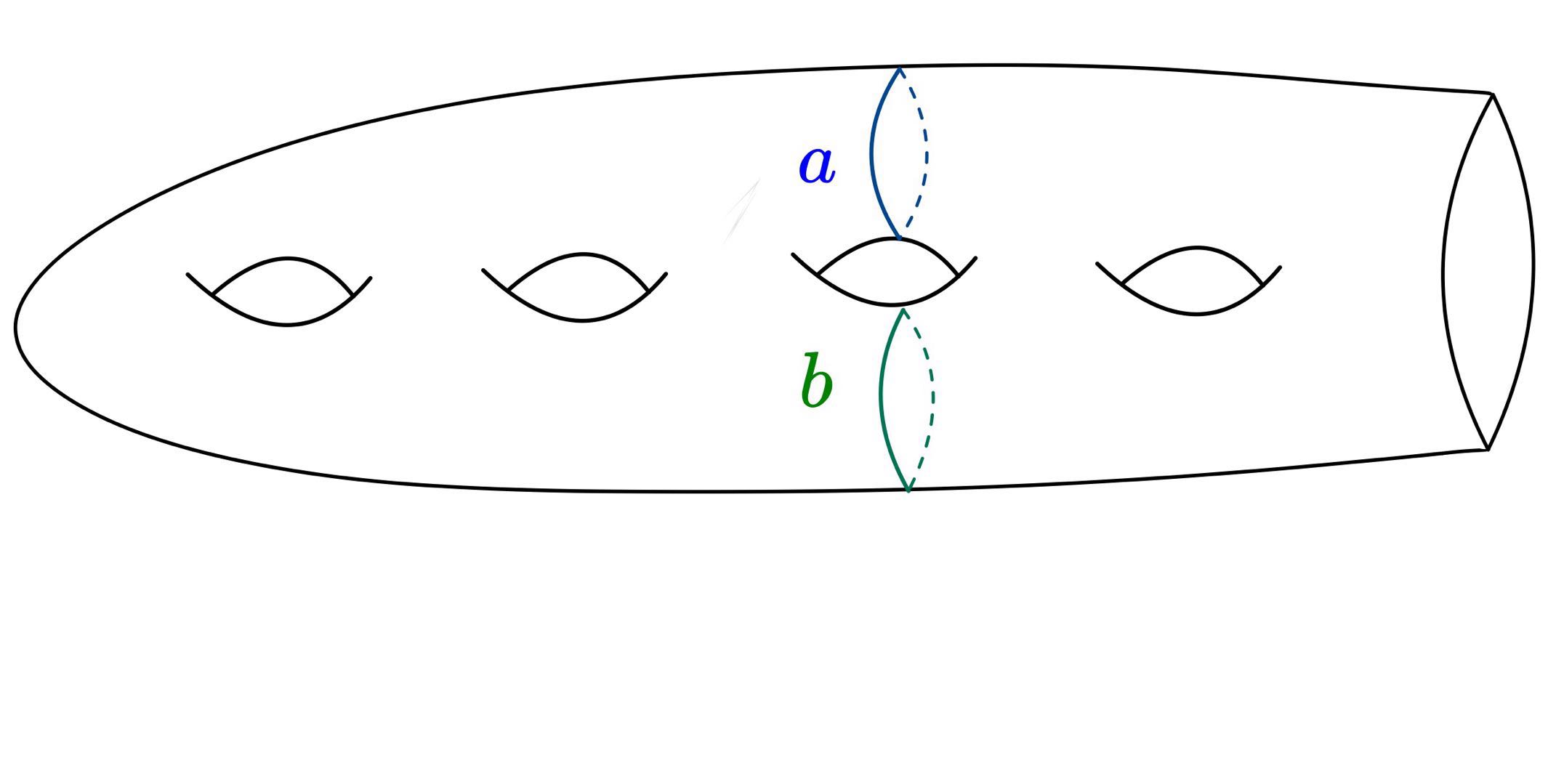}
        \vspace{-35pt}
    \caption{$T_aT_b^{-1}$ is a genus 2 bounding pair map.}
    \label{fig:placeholder}
\end{figure}
        \vspace{-7pt}
        
Let $\W_g^1(n)$ denote the normal subgroup generated by the genus $n$ bounding pair maps in the mapping class group $\M_g^1$ of the once-bordered surface $S_g^1$ of genus $g$. Let $\ch:\I_g^1\xrightarrow{\ch} H_1(S_g^1;\Z)$ denote the Chillingworth homomorphism on the Torelli group $\I_g^1$. Define a \emph{Chillingworth congruence subgroup} to be
$$\Ch_g^1[2n]:=\{f\in\I_g^1\ : \ \ch(f)=0\pmod{2n}\}.$$

\begin{theorem}[\textbf{First description of $\mathbf{\W_g^1(n)}$}]
\label{thm:W is Ch M commutator}

    When $1\le n\le g-2$, we have that
    $$\W_g^1(n)=[\Ch_g^1[2n],\M_g^1].$$
    
\end{theorem}
As corollaries, we have that $\W_g^1(n)$ is rich enough to contain, for example, the commutator of any separating twist with any mapping class. On the other hand, in contrast to Johnson's result that $\W_g^1(1)=\I_g^1$, we have 
$\W_g^1(n)\subseteq \Ch_g^1[2n]\varsubsetneqq \I_g^1$ when $n>1$. Our next theorem determines the index of $\W_g^1(n)\subseteq\Ch_g^1[2n]$. 

\begin{theorem}[\textbf{Second description of $\mathbf{\W_g^1(n)}$}]
\label{thm: W is kerd}
Let $d_{4n}: \Ch_g^1[2n]\to\Z/4n\Z$ denote the Casson--Morita map modulo $4n$. 
When $1\le n\le g-2$, the map
$d_{4n}$ is an $\M_g^1$-invariant homomorphism such that    
\[\W_g^1(n)=\ker(d_{4n}).\]
In these cases, we have that 

\[
\frac{\Ch_g^1[2n]}{\W_g^1(n)} \cong\begin{cases}
			\Z/n\Z, & \text{if $n$ is odd}\\
            \Z/\frac{n}{2}\Z, & \text{if $n$ is even}.
		 \end{cases}
\]
\end{theorem}

Let us illustrate Theorem \ref{thm: W is kerd} with a few examples. Let $T_h$ denote a genus $h$ separating twist. Morita \cite[Theorem 5.3]{Morita2} showed that $d(T_h)=4h(h-1)$.  Hence, Theorem \ref{thm: W is kerd} implies that $T_1\in\W_g^1(n)$, and that $T_2\in\Ch_g^1[2n]$ but $T_2\notin\W_g^1(n)$ when $n>2$. On the other hand, the commutator of $T_2$ with any mapping class is in $\W_g^1(n)$. 

Theorem \ref{thm: W is kerd} is similar in spirit to another theorem of Johnson which says that the group $\J_g^1$ generated by all separating twists is equal to the kernel of the Johnson homomorphism. Both Theorem \ref{thm: W is kerd} and Johnson's result characterize a geometrically defined subgroup of the Torelli group as the kernel of an algebraic invariant. As an application of Theorem~\ref{thm: W is kerd}, we prove how much of $\J_g^1$ the subgroup $\W_g^1(n)$ contains.

\begin{corollary}
    \label{cor:W vs K}
    $\W_g^1(n)$ contains $\J_g^1$ if and only if $n=1,2$. More generally, for $1\le n\le g-2$, we have that
    $$\frac{\J_g^1}{\W_g^1(n)\cap \J_g^1} \cong\begin{cases}
			\Z/n\Z, & \text{if $n$ is odd}\\
            \Z/\frac{n}{2}\Z, & \text{if $n$ is even}.
		 \end{cases}$$
\end{corollary}

Note that the hypothesis $n\le g-2$ for these results is necessary. For example, Theorem \ref{thm: W is kerd} implies that $\W_g^1(n)$ is always a finite index subgroup of $\I_g^1$ when $1\le n\le g-2$. In contrast, $\W_g^1(g-1)$ is equal to the point-pushing subgroup which is of infinite index in $\I_g^1$ (see \cite[Section 3]{johnson1983structure}).

As corollaries, we compute the images of $\W_g^1(n)$, $\Ch_g^1[2n]$, and $\I_g^1$ under $d$ in Section \ref{sec: integral d}. We also deduce the following homological results:
\begin{corollary}
        \label{cor:homological results}
    When $1\le n\le g-2$, we have that
    \begin{enumerate}
        \item    $H_1(\Ch_g^1[2n];\Z)_{\M_g^1}\cong
        \begin{cases}
			\Z/n\Z, & \text{if $n$ is odd}\\
            \Z/\frac{n}{2}\Z, & \text{if $n$ is even}.
		 \end{cases}
         $
         \item $[\W_g^1(n),\M_g^1]=\W_g^1(n)$, which implies that  $H_1(\W_g^1(n);\Z)_{\M_g^1}=0$. 
    \end{enumerate}
\end{corollary}

In Section \ref{sec:closed and punctured} we obtain similar results for once-punctured and closed surfaces. 

\begin{corollary}
    \label{cor: punctured and closed case in the intro}
    Suppose that $1 \leq n\leq g-2$. 
    \begin{enumerate}
        \item (once-punctured surfaces) For $\ell:=\gcd(g(g-1),n)$, 
\[
\W_{g,1}(n)=\ker(d_{4\ell}: \Ch_{g,1}[2n]\to\Z/4\ell\Z).
\]
In this case, we have that  
    \[
\frac{\Ch_{g,1}[2n]}{\W_{g,1}(n)}\cong\frac{\J_{g,1}}{\W_{g,1}(n)\cap \J_{g,1}} \cong\begin{cases}
			\Z/\ell\Z, & \text{if $n$ is odd}\\
            \Z/\frac{\ell}{2}\Z, & \text{if $n$ is even}.
		 \end{cases}
\]    
In particular, $\W_{g,1}(n)=\I_{g,1}$ if and only if $n=1$; $\W_{g,1}(n)$ contains $\J_{g,1}$ if and only if $\ell\in \{1,2\}$. 
\item (closed surfaces) For $m:=\gcd(g-1,n)$, 
\[
\W_{g}(n)=\ker(d_{4m}: \Ch_{g}[2n]\to\Z/4m\Z).
\]
In this case, we have that  $\W_g(n)=\W_g(m)$ and that
    \[
\frac{\Ch_{g}[2n]}{\W_{g}(n)}\cong\frac{\J_{g}}{\W_{g}(n)\cap \J_{g}} \cong\begin{cases}
			\Z/m\Z, & \text{if $m$ is odd}\\
            \Z/\frac{m}{2}\Z, & \text{if $m$ is even}.
		 \end{cases}
\]    
In particular, $\W_g(n)=\I_g$ if and only if $n$ and $g-1$ are relatively prime; $\W_g(n)$ contains $\J_g$ if and only if $m\in\{1,2\}$.
\end{enumerate}
\end{corollary}

\p{Outline}
In Section \ref{sec:image johnson}, we compute the image of $\W_g^1(n)$ under the Johnson homomorphism. Then in Section \ref{sec:septwists}, we show that various separating twists and products thereof are elements in $\W_g^1(n)$ using the lantern relation in different ways. A key step is to show that a certain curve complex is connected using a lemma of Putman. In Section \ref{sec:proof of main}, we discuss the Casson--Morita map $d$ and prove the main results. In Section \ref{sec:closed and punctured}, we prove similar results for once-punctured and closed surfaces.

\p{Acknowledgments} The authors would like to thank Dan Margalit and Marissa Loving for the conversations that set this work in motion, and to thank Dan Margalit for comments on the article. The first author is supported by the National Science Foundation under Grant No. DMS-2203178 and by a Sloan Fellowship. The second author is supported by the National Natural Science Foundation of China under the Young Scientists Fund No. 12101349. The third author is supported by the National Science Foundation under Grant No. DGE-1650044 and Grant No. DMS-2002187. We thank an anonymous referee for their helpful comments on this article. In an earlier version of this paper, we made a mistake in the computation of the image of $\W_g^1(n)$ under the Johnson homomorphism. We thank the anonymous referee for pointing this out.

\section{Image of $\W_g^1(n)$ under the Johnson homomorphism}
\label{sec:image johnson}
\subsection{Chillingworth and Johnson  homomorphisms}
For now, let us focus on the  oriented surface $S_g^1$ of genus $g$ with 1 boundary component. Set $H:=H_1 (S_g^1;\Z)$. The Torelli group $\I_g^1$ is the subgroup of the mapping class group $\M_g^1$ that acts trivially on $H$. Chillingworth \cite{chillingworth1972winding} studied a crossed-homomorphism 
$\ch: \M_g^1\to H$
defined by the action of $\M_g^1$ on vector fields on $S_g^1$ up to homotopy.
Johnson \cite{JohnsonAbelian} proved that the restriction of the Chillingworth homomorphism $\ch$ to $\I_g^1$ factors through a surjective homomorphism 
\[
\tau_g^1: \I_g^1\to \wedge^3 H.
\]
Specifically, Johnson \cite[Theorem 2]{JohnsonAbelian} proved that 
\begin{equation}
    \label{eq:ch=C tau}
    e(f)=C(\tau_g^1(f))\ \ \ \ \ \ \ \ \ \forall f\in \I_g^1
\end{equation}
where $C: \wedge^3 H \to H$ is the contraction homomorphism  defined by 
\begin{equation}
    \label{eq: C def}
    C(x \wedge y \wedge z) =2\Big[(x\cdot y) z+ (y\cdot z) x+ (z\cdot x) y\Big]
\end{equation}
where $x\cdot y$ denotes the intersection pairing of $x$ and $y$. For convenience, we will use (\ref{eq:ch=C tau}) as the definition of the {Chillingworth homomorphism} on $\I_g^1$. The \emph{Chillingworth subgroup} is $\Ch_g^1:= \text{Ker}(\ch)\le \I_g^1$.

\subsection{Calculating $\tau_g^1(\W_g^1(n))$}
We will focus on the bordered surface $S_g^1$, and therefore we will simply use $\M, \I,\tau,\W(n)$ to denote $\M_g^1, \I_g^1,\tau_g^1, \W_g^1(n)$, respectively.

Let $\{a_1,b_1,\dots, a_g,b_g\}$ denote a symplectic basis for $H$ under the intersection pairing.  Let $BP_n$ denote the standard bounding pair map of genus $n$ with respect to this basis. A straightforward calculation (see \emph{e.g.} \cite[Section 6.6.2]{primer} for details) gives that
\[
\tau(BP_n)=\sum_{i=1}^n (a_i \wedge b_i) \wedge b_{n+1}.
\]

Applying Johnson's contraction homomorphism $C$ defined in (\ref{eq: C def}), we obtain
\begin{equation}
    \label{ch of BPn}
    \ch(BP_n)=C(\tau(BP_n))=2n (b_{n+1})\in 2nH.
\end{equation}

This implies that 
\begin{equation}
    \label{tau(W_n)<U2n}
    \tau(\W(n))\le U_{2n}:= C^{-1}(2nH).
\end{equation}
Our goal in this section is to prove that the containment (\ref{tau(W_n)<U2n}) is actually an equality.

Before we proceed, let us remark that (\ref{tau(W_n)<U2n}) implies that $\W(n)$ is a proper subgroup of $\I$ when $n>1$ because
\begin{equation}
    \label{eq: W in Ch}
    \W(n)\le \Ch_g^1[2n]=e^{-1}(2nH).
\end{equation}
However, we will prove later that this containment is generally not an equality because certain elements in the Johnson kernel, such as any genus 2 separating twists, are not in general in $\W(n)$.

\begin{proposition}
\label{generateU}
The subgroup $U:=\ker(C)$ is generated by elements of the following form: 
\begin{enumerate}
\item $a_i\wedge a_j\wedge a_k$ for distinct $i,j,k\in \{1,...,g\}$,
\item $b_i\wedge b_j\wedge b_k$ for distinct $i,j,k\in \{1,...,g\}$,
\item $a_i\wedge a_j\wedge b_k$ for distinct $i,j,k\in \{1,...,g\}$,
\item  $a_i\wedge b_j\wedge b_k$ for distinct $i,j,k\in \{1,...,g\}$,
\item $(a_i\wedge b_i-a_j\wedge b_j)\wedge a_k$ for distinct $i,j,k\in \{1,...,g\}$,
\item $(a_i\wedge b_i-a_j\wedge b_j)\wedge b_k$ for distinct $i,j,k\in \{1,...,g\}$.
\end{enumerate}
\end{proposition}
\begin{proof}
$\wedge^3 H$ is generated by the elements in (1)-(6) above together with elements of the following form:
\begin{enumerate}
\item[\emph{(7)}] $a_1\wedge b_1\wedge a_{k}$   for $k\in \{2,...,g\}$ and $a_2\wedge b_2\wedge a_1$,
\item[\emph{(8)}] $a_1\wedge b_1\wedge b_{k}$   for $k\in \{2,...,g\}$ and $a_2\wedge b_2\wedge b_1$.
\end{enumerate}
The contraction map $C$ takes generators in (1)-(6) to $0$, and the $2g$  generators in (7) and (8) into a standard basis for $2H\subseteq H$. Hence,  $U=\ker C$ is spanned by elements in (1)-(6).
\end{proof}

\begin{proposition}
\label{U2n generators}
For any $n \ge 1$, the subgroup
$U_{2n}:=C^{-1}(2nH)$
is generated by $U$ along with the image of $\tau(BP_n)$ under the natural action of $\Sp$, where $BP_n$ denotes the standard bounding pair map of genus $n$.
\end{proposition}

\begin{proof}
    Our calculation (\ref{ch of BPn}) implies that the subgroup generated by the image of $C(\tau(BP_n))$ under the $\Sp$-action is
    $$\Z \Sp\{ C(\tau(BP_n))\}=2nH$$
    because $\Z \Sp\{ b_{n+1}\}=H$. In particular, we have that $C(U_{2n})= 2nH$.  The proposition follows from the following short exact sequence
\[0\to U\to U_{2n}\xrightarrow{C}  2nH\to 0. \qedhere\]
    \end{proof}
We next calculate the image of $\W(n)$ under the Johnson homomorphism.
\begin{proposition}
\label{tau(W_n)=U2n}
For $1\le n\le g-2$, we have that $\tau(\W(n))=\tau(\Ch_g^1[2n])=U_{2n}$.
\end{proposition}
\begin{proof}
We know that 
$\tau(\W(n))\le\tau(\Ch_g^1[2n])\le U_{2n}$
by (\ref{tau(W_n)<U2n}). Therefore, it suffices to prove that $U_{2n}\le\tau(\W(n))$. By Proposition \ref{U2n generators}, it remains to show that $U\le \tau(\W(n))$.  We will check that each of the generators of $U$ in Proposition \ref{generateU} is in $\tau(\W(n))$.

Let $\phi$ be the factor mix element of $\Sp$ defined by the following map on generators of $H$:
\begin{align*}
a_1 \mapsto a_1 - b_{n+2}\\
a_{n+2} \mapsto a_{n+2} - b_1
\end{align*}
and fixing all other basis elements. This $\phi$ exists because $n\le g-2$. Again, let $BP_n$ denote the standard bounding pair map of genus $n$.
\[
\tau(BP_n)=\sum_{i=1}^n (a_i \wedge b_i) \wedge b_{n+1}\in \tau (\W(n))
\]

\[
\phi(\tau(BP_n))=(-b_{n+2} \wedge b_1 \wedge b_{n+1}) + \sum_{i=1}^n ( a_i \wedge b_i) \wedge b_{n+1} \in  \tau (\W(n)).
\]
Hence, we have
\[
\phi(\tau(BP_n))-\tau(BP_n)=(-b_{n+2} \wedge b_1 \wedge b_{n+1}) \in  \tau (\W(n)).
\]
By the action of $\Sp$, we obtain elements of the forms (1)-(4) in Proposition \ref{generateU}.

Similarly, let $\psi$ be the factor swap element of $\Sp$ given by the following map on generators of $H$:
\begin{align*}
a_1 \mapsto a_{n+2} & & a_{n+2} \mapsto a_1\\
b_1 \mapsto b_{n+2} & & b_{n+2} \mapsto b_1
\end{align*}
and fixing all other basis elements. This $\psi$ exists by the hypothesis that $n+2\le g$. 
\[
\psi(\tau(BP_n))-\tau(BP_n)=(a_1\wedge b_1-a_{n+2}\wedge b_{n+2})\wedge  b_{n+1} \in \tau (\W(n)).
\]
By the action of $\Sp$, we obtain elements of the form (5) and (6) in Proposition \ref{generateU}. 
\end{proof}

\section{Separating twists in $\W(n)$}
\label{sec:septwists}

In this section, we consider the surface $S_g^1$ and will again suppress the underlying surface from notations.  Let $\T(m)$ be the normal closure of a genus $m$ separating twist in $\M$. We will show that $\W(n)$ contains $\T(m)$ for $m\in\{1,n,n+1\}$. We will also show that $\W(n)$ contains all elements of the form $T_bT_c^{-1}$ where $b,c$ are any separating curves of genus $2$. Several proofs in this section are similar to arguments made by Johnson in showing that $\langle \T(1), \T(2) \rangle$ is equal to the subgroup generated by all separating twists and that $\W(1)=\I$ \cite[Theorems~1~and~2]{JohnsonLantern}, using lantern relations.

\begin{proposition}
\label{prop:BPtosep}
For the surface $S_g^1$, we have that
\begin{enumerate}
    \item $\T(n) \le \W(n)$ for $1\leq n \le g-2$
    \item $\T(n+1) \leq \W(n)$  for $1\leq n \le g-1$
    \item $\T(1) \leq \W(n)$  for $1\leq n \le g-2$
\end{enumerate}
\end{proposition}
\begin{proof}
The conditions on $g$ and $n$ in \textit{(1)} and \textit{(2)} guarantee that $S$ has a lantern in the two configurations depicted in Figure~\ref{fig:BPtosep}, respectively. The two configurations differ only in the location of the boundary component. 
\begin{figure}[b]
\centering
{\labellist\small\hair 2.5pt

        \pinlabel $\boxed{n}$ by 0 0 at 30 170
         \pinlabel $\ast$ by 0 0 at 280 170
        \pinlabel $\boxed{g-n-2}$ by 0 0 at 263 140
 \pinlabel $x$ by 0 0 at 140 200
  \pinlabel $y$ by 0 0 at 140 130
   \pinlabel $z$ by 0 0 at 170 170
      \pinlabel $c$ by 0 0 at 212 170
              \pinlabel $d$ by 0 0 at 62 170
               \pinlabel $a$ by 0 0 at 130 250
  \pinlabel $b$ by 0 0 at 130 70
  
            \pinlabel $\Big\uparrow$ by 0 0 at 370 105
          \pinlabel $\boxed{g-n-1}$ by 0 0 at 370 70
         \pinlabel $\ast$ by 0 0 at 360 170
 \pinlabel $x$ by 0 0 at 480 200
  \pinlabel $y$ by 0 0 at 480 130
   \pinlabel $z$ by 0 0 at 510 170
      \pinlabel $c$ by 0 0 at 552 170
      \pinlabel $\boxed{n-1}$ by 0 0 at 600 140
              \pinlabel $d$ by 0 0 at 402 170
               \pinlabel $a$ by 0 0 at 470 250
  \pinlabel $b$ by 0 0 at 470 70
  
        \endlabellist}
\includegraphics[width=.48\textwidth]{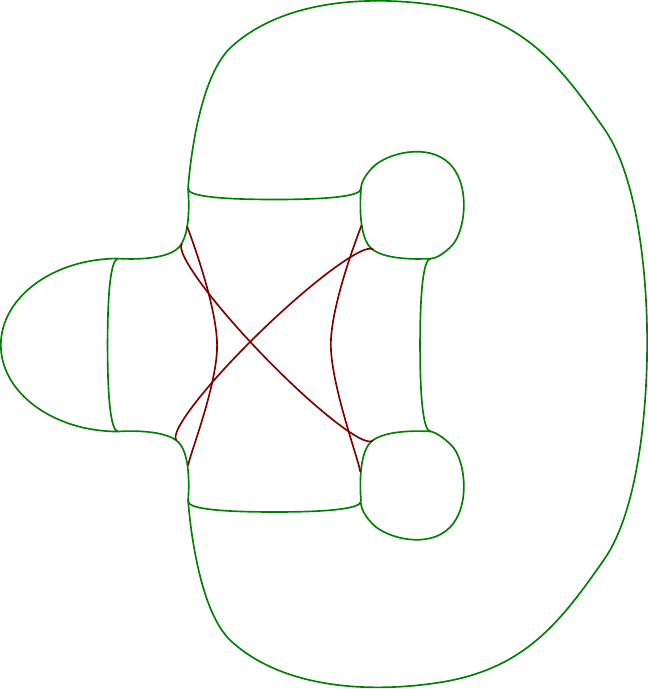} \ \ \
\includegraphics[width=.48\textwidth]{BPtosepclean}
\caption{Two lantern configurations used to show that $\T(n)$ and $\T(n+1)$ respectively are subgroups of $\W(n)$. The asterisks indicate the boundary components, and the boxed values in the subsurfaces indicate their genus. }
\label{fig:BPtosep}
\end{figure}

In both figures, we have the lantern relation:
\[
T_aT_bT_cT_d = T_xT_yT_z
\]
or 
\[
T_d=(T_xT_a^{-1})(T_yT_b^{-1})(T_zT_c^{-1}).
\]
In the first configuration, $T_d$ is a separating twist of genus $n$, while each of the three terms on the right is a bounding pair map of genus $n$. We therefore conclude that $\T(n) \leq \W(n)$. In the second configuration, we have that $T_d$ is a separating twist of genus $n+1$, while each of the three terms on the right is a bounding pair map of genus $n$. We therefore conclude that $\T(n+1) \leq \W(n)$.

To prove $\T(1)\leq\W(n)$, we take a lantern in $S$ as depicted on the left in Figure~\ref{fig:BPtoJ}. 
\begin{figure}[b]
\centering
{\labellist\small\hair 2.5pt

        \pinlabel $\boxed{n}$ by 0 0 at 30 130
         \pinlabel $\ast$ by 0 0 at 130 10
 \pinlabel $x$ by 0 0 at 140 165
  \pinlabel $y$ by 0 0 at 140 95
   \pinlabel $z$ by 0 0 at 170 130
      \pinlabel $c$ by 0 0 at 212 130
              \pinlabel $d$ by 0 0 at 62 130
               \pinlabel $a$ by 0 0 at 130 210
  \pinlabel $b$ by 0 0 at 130 40
    \pinlabel $\boxed{g-n-1}$ by 0 0 at 242 20
                \pinlabel $\longleftarrow$ by 0 0 at 184 20
                
          \endlabellist}
\includegraphics[width=.49\textwidth]{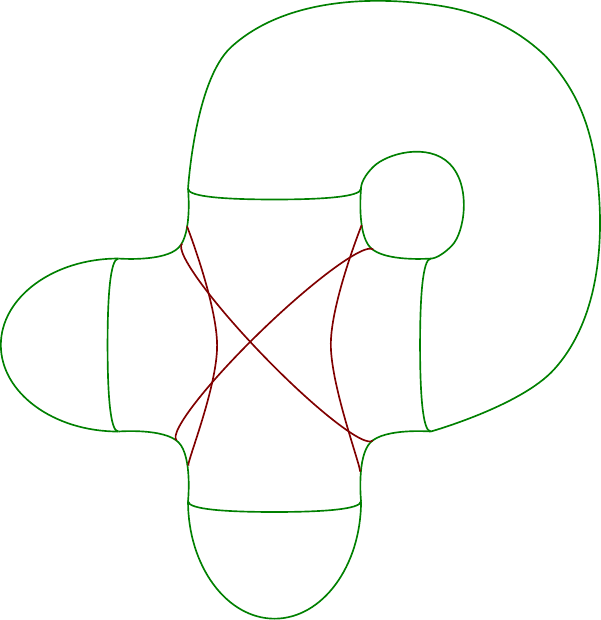}
{\labellist\small\hair 2.5pt

        \pinlabel $\boxed{n}$ by 0 0 at 28 130
                \pinlabel $\boxed{m}$ by 0 0 at 240 230
         \pinlabel $\ast$ by 0 0 at 130 10
 \pinlabel $x$ by 0 0 at 140 165
  \pinlabel $y$ by 0 0 at 140 95
   \pinlabel $z$ by 0 0 at 170 130
      \pinlabel $c$ by 0 0 at 212 130
              \pinlabel $d$ by 0 0 at 62 130
               \pinlabel $a$ by 0 0 at 130 210
  \pinlabel $b$ by 0 0 at 130 40
    \pinlabel $\boxed{g-n-m}$ by 0 0 at 242 20
                \pinlabel $\longleftarrow$ by 0 0 at 184 20
                
          \endlabellist}
\includegraphics[width=.49\textwidth]{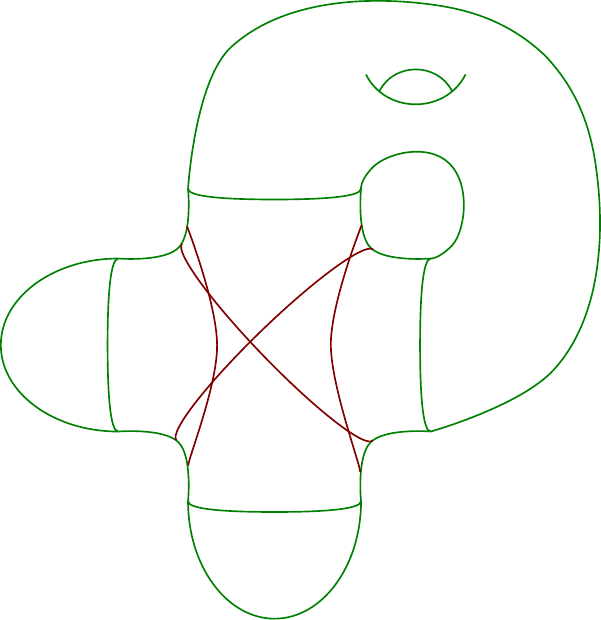}
\caption{Left: the lantern configuration used to show that $\T(1)$ is a subgroup of $\W(n)$. Right: the lantern configuration used to prove Proposition \ref{prop:nested}. The asterisks indicate the boundary components, and the boxed values in the subsurfaces indicate their genus. }
\label{fig:BPtoJ}
\end{figure}
The separating curves $b$, $d$, and $y$ cut off subsurfaces of genera $n+1$, $n$, and $1$, respectively. By the lantern relation, we have
\[
T_y=(T_aT_x^{-1})(T_cT_z^{-1})T_bT_d.
\]
We have that $T_y$ is a separating twist of genus $1$. Each of the first two terms on the right is a bounding pair map of genus $n$ and so these lie in $\W(n)$. The other two terms are separating twists of genera $n+1$ and $n$, respectively; these lie in $\W(n)$ by \textit{(1)} and \textit{(2)}. We therefore conclude that $\T(1)\leq \W(n)$.
\end{proof}

Using a similar argument, we obtain the following.
\begin{proposition}
\label{prop:nested}
For the surface $S_g^1$,  when $1\leq n \le g-2$ and $m+n \leq g$, we have that $T_yT_b^{-1} \in \W(n)$ where $y$ and $b$ are nested separating curves of genus $m$ and $m+n$, respectively.
\end{proposition}

\begin{proof}
We may take a lantern in $S$ as depicted on the right in Figure~\ref{fig:BPtoJ}. The separating curves $b$, $d$, and $y$ cut off subsurfaces of genera $m+n$, $n$, and $m$, respectively. By the lantern relation, we have
\[
T_yT_b^{-1}=(T_aT_x^{-1})(T_cT_z^{-1})T_d.
\]
We have that the separating curves $y$ and $b$ are nested, that $T_y$ is a separating twist of genus $m$, and that $T_b^{-1}$ is a separating twist of genus $m+n$. Each of the first two terms on the right is a bounding pair map of genus $n$ and so these lie in $\W(n)$. The separating twist $T_d$ of genus $n$ lies in $\W(n)$ by Proposition~\ref{prop:BPtosep}. We therefore conclude that $T_yT_b^{-1} \in \W(n)$. By the classification of surfaces, all products of nested separating twists of genus $m$ and $m+n$ are conjugate in $\M$.\end{proof}

\begin{corollary}
\label{genus2}
For the surface $S_g^1$, when $g \geq 4$ and $1\le n \le g-2$, we have that $T_bT_c^{-1}\in \W(n)$ where $b,c$ are disjoint separating curves of genus $2$.
\end{corollary}
\begin{proof}
Under the hypothesis, $S$ contains three disjoint separating curves $b,c,y$ of genera $2,2,n+2$ respectively such that $b,y$ and $c,y$ are both nested. Then we have that $T_bT_y^{-1}\in \W(n)$ and $T_cT_y^{-1}\in \W(n)$. Thus we have that $T_bT_c^{-1}\in \W(n)$.
\end{proof}

Next, we will strengthen Corollary \ref{genus2} by removing the hypothesis that $b$ and $c$ are disjoint, provided that $g\ge5$.
\begin{proposition}
    \label{prop:difference of T2 in W}
    For the surface $S_g^1$, when $g \geq 5$ and $1\le n \le g-2$, we have that $T_bT_c^{-1}\in \W(n)$ where $b,c$ are any separating curves of genus $2$.
\end{proposition}

Before proving this proposition, let us introduce the following lemma due to Putman \cite[Lemma 2.1]{putman2008note}.

\begin{lemma}[Putman's connectivity lemma]\label{Putman}
Consider a group $G$ acting upon a simplicial complex $X$ where $X_n$ denotes its $n$th skeleton. Fix a basepoint $v \in X_0$ and a set $S$ of generators for $G$. Assume the following hold.
\begin{enumerate}
\item For all $v' \in X_0$, the orbit $Gv$ intersects the connected component of $X$ containing $v'$;
\item  For all $s \in S^{\pm 1}$, there is some path $P(s)$ in $X$ from $v$ to $sv$.
\end{enumerate}
Then X is connected.
\end{lemma}

Let $\X_{g,\text{sep}}^{1}(2)$ be the curve complex consisting of separating curves of genus $2$ as vertices. Two curves are connected by an edge if and only if they are disjoint. In particular, $\X_{g,\text{sep}}^{1}(2)$ is the full subcomplex of the curve complex $\X_g^1$ spanned by separating curves bounding a genus 2 subsurface in $S_g^1$. Note that the mapping class group $\M_g^1$ acts naturally on $\X_{g,\text{sep}}^{1}(2)$.
\begin{proposition}
\label{con}
For $g\ge 5$, the curve complex $\X_{g,\text{sep}}^{1}(2)$ is connected. 
\end{proposition}
\begin{proof}
The group $\M_g^1$ acts on $\X_{g,\text{sep}}^{1}(2)$ transitively.  We have Lickorish's generating set for $\M_g^1$ as in Figure \ref{Lick}. 
\begin{figure}[b]
\centering
\includegraphics[width=.8\textwidth]{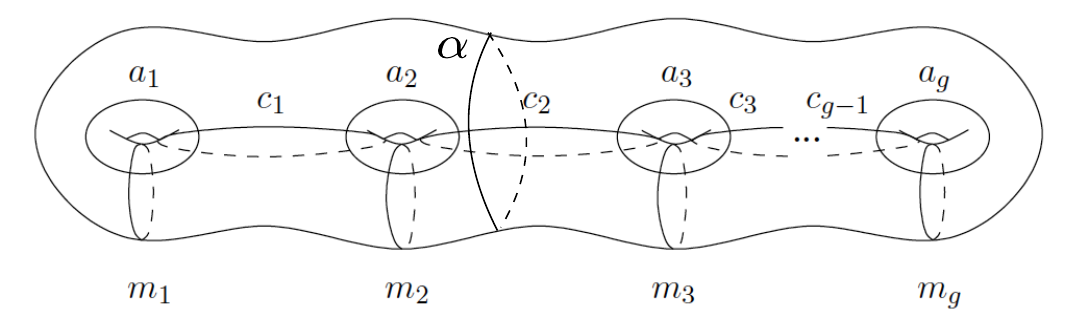}
\caption{Lickorish's generating set from Figure 5.7 in \cite{primer} with the curve $\alpha$ added.}
\label{Lick}
\end{figure}
Let $v:= \alpha$ be the curve in Figure \ref{Lick}. If $g\ge 5$, there exists a simple closed  separating curve $\beta$ of genus $2$ distict and disjoint from both $\alpha$ and $c_2$. Since the only curve in $\{a_1,...,m_g\}$ that intersects $\alpha$ is $c_2$, and $T_{c_2}(\alpha)$ is disjoint from $\beta$, we  have both (1) and (2) of Lemma \ref{Putman}.
\end{proof}
The above connectivity result gives us Proposition \ref{prop:difference of T2 in W}.

\begin{proof}[Proof of Proposition \ref{prop:difference of T2 in W}]
   When $g\ge 5$, by Proposition \ref{con}, we know that $b$ and $c$ are connected by a path $(a_0,\dots,a_k)$ in $\X_{g,\text{sep}}^{1}(2)$ such that $a_0=b$ and $a_n=c$. Then we can write 
    \[
T_bT_{c}^{-1} =(T_{a_0}T_{a_{1}}^{-1})\dots (T_{a_{n-1}}T_{a_{n}}^{-1})
    \]
Since each term on the right is a product of disjoint genus 2 separating twists, and these lie in $\W(n)$ by Corollary~\ref{genus2}, we conclude that $T_bT_{c}^{-1}\in \W(n)$.
\end{proof}

\section{Proving the main results}
\label{sec:proof of main}

In this section, we again will focus on the bordered surface $S_g^1$ and sometimes suppress the underlying surface from notation. 

\subsection{The Casson--Morita map and its restriction to $\Ch_g^1[2n]$}
\label{sec:Def of d}
The most important tool in the proofs of the main theorems is the Casson--Morita map $d:\M_g^1\to \Z$. We  will first briefly recall the construction of this map and its properties following Morita's paper \cite{Morita2}. Next, we will prove several key properties of this map when restricted to $\Ch_g^1[2n]$.

Morita constructed a crossed homomorphism $k: \M_g^1 \to H:=H_1(S_g^1;\Z)$ in the following way. Define a map $\epsilon:\pi_1(S_g^1)\to\Z$ by $\epsilon(\gamma)= \sum_{i=1}^{r-1} [\gamma_1]\cdot [\gamma_{i+1}...\gamma_r]$ where $\gamma=\gamma_1...\gamma_r$ is the unique expression of an element $\gamma$ in terms of the standard generators $\alpha_i^{\pm1}$ and $\beta_i^{\pm1}$. Any $\phi\in\M_g^1$ gives a homomorphism $k^\vee(\phi):H\to\Z$ defined by $k^\vee(\phi)(\gamma):=\epsilon(\phi(\gamma))-\epsilon(\gamma)$. Then we define $k(\phi)$ by the equation $k^\vee(\phi)(u)=k(\phi)\cdot u$ for any $u\in H$ where $\cdot$ denotes the intersection form. It is straightforward to check that 
\begin{equation}
    \label{eq: k properties}
    k(\phi\psi)=\psi^{-1}(k(\phi))+k(\psi)\ \ \ \ \text{ for any $\phi,\psi\in\M_g^1$}. 
\end{equation}
In other words, $k$ is a crossed homomorphism or equivalently a 1-cocyle. In fact, Morita showed that $k$ represents a generator for $H^1(\M_g^1;H)\cong \Z$. Furthermore, define a 2-cocyle by $c(\phi,\psi):=k(\phi)\cdot k(\psi^{-1})$. Let $\sigma$ denote \emph{Meyer's signature cocycle} $\sigma:\Sp\times\Sp\to\Z$ pulled back to $\M_g^1.$ Morita showed that $c+3\sigma$ represents the zero cohomology class in $H^2(\M_g^1;\Z)$. Consequently, there exists a function $d:\M_g^1\to \Z$ such that
\begin{equation}
    \label{eq:d=}
    \delta d=c+3\sigma.
\end{equation}
This function $d$ is called the \emph{Casson--Morita map} because the restriction of $d$ to $\J_g^1$ is related to the Casson invariant of homology 3-spheres \cite[Theorem 6.1]{Morita2}. Such $d$ is unique when $g\ge 3$ since $H^1(\M_g^1;\Z)=0$ in this case.

\begin{proposition}[Morita, Proposition 5.1 and Theorem 5.3 in \cite{Morita2}]
\label{d map properties}
For any two elements $\phi,\psi \in \M_g^1$, 
\begin{enumerate}
\item 
$d ( \phi \psi ) = d ( \phi ) + d ( \psi ) +k( \phi ) \cdot \psi _ { * } k( \psi ) - 3 \sigma ( \phi , \psi )$;
\item
$d ( \phi ^ { - 1 } ) = - d ( \phi )$;
\item
$d ( \phi \psi \phi ^ { -1 } ) = d ( \psi ) + k( \phi ) \cdot( \psi _ { * } k( \psi ) + k ( \phi \psi ) )$
\end{enumerate}
Moreover, if $T_h$ denotes a separating twist of genus $h$, then $d(T_h)=4h(h-1).$
\end{proposition}

Morita \cite[Theorem 6.1]{Jacobi2} showed that the restriction of  $k$ to $\I_g^1$ is equal to the Chillingworth homomorphism $\ch=C\circ\tau$ up to a sign. 
 Hence,  Proposition \ref{d map properties} implies that $d$ is an $\M_g^1$-invariant homomorphism over $\Ch_g^1=\ker \ch$. However, in general $d$ is not a homomorphism even when restricted to $\I_g^1$.

\begin{lemma}
\label{d4n invariant homo}
For $m\in\Z$, define $d_{m}:\Ch_g^1[2n]\to \Z/m\Z$ as the composition:
$$\Ch_g^1[2n]\xrightarrow{d|_{\Ch_g^1[2n]}}\Z\to \Z/m\Z.$$
Then we have that
\begin{enumerate}
\item $d_m$ is a homomorphism if $m=(2n)^2$;
\item  $d_{4n}$ is a $\M_g^1$-invariant homomorphism. 
\end{enumerate}
\end{lemma}
\begin{proof}
Since $k$ coincides with the Chillingworth homomorphism $\ch$ up to a sign, $k(\Ch_g^1[2n])=\ch(\Ch_g^1[2n])\subseteq 2nH$.
For any $\phi,\psi\in \Ch_g^1[2n]$, Proposition \ref{d map properties} gives
\[
d ( \phi \psi ) = d ( \phi ) + d ( \psi ) + k ( \phi ) \cdot k( \psi )-0.
\]
Notice that since $\phi,\psi$ are in the Torelli group, we have that $\psi_*=id$ and that $\sigma(\phi,\psi)=0$. Since $k(\phi),k(\psi)\in 2nH$, the difference $d ( \phi \psi ) - d ( \phi ) - d ( \psi )$ is a multiple of $(2n)^2$.

For any $\phi\in \M_g^1$ and $\psi\in \Ch_{g,1}[2n]$, Proposition \ref{d map properties} gives
\begin{equation}
\begin{aligned}
d ( \phi \psi \phi ^ { -1 } ) &= d ( \psi ) + k( \phi ) \cdot ( k ( \psi ) + k ( \phi \psi ) ) \\
&= d ( \psi ) + k ( \phi ) \cdot( k( \psi ) +k ( \phi )+k(\psi) ) )
 \\
&=d(\psi)+2k(\phi)\cdot k(\psi)\\
\end{aligned}
\end{equation}
where $k(\phi)\cdot k(\psi)\in 2n\Z$. We thus know that $d_{4n}$ is $\M_g^1$-invariant.
\end{proof}

\subsection{Proof of Theorem \ref{thm:W is Ch M commutator}}
Theorem \ref{thm:W is Ch M commutator} will follow from Theorem \ref{main1} below.
\begin{lemma} 
\label{Kerd}
For $1\le n\le g-2$, we have $\ker(d|_{\J_g^1})\le \W(n)$.\end{lemma}
\begin{proof}
Faes \cite[Remark 2.15]{Fae} showed that $\ker(d|_{\J_g^1})$ is normally generated by $[\J_g^1,\M_g^1]$ and a genus $1$ separating twist. By Johnson \cite{JohnsonAbelian}, we know that $\J_g^1$ is generated by separating twists of genera 1 and 2. Hence, $\ker(d|_{\J_g^1})$ is normally generated by a genus 1 separating twist and $[T_a,f]$ for all $f\in \M_g^1$ and genus 2 separating curves $a$. 

Proposition \ref{prop:BPtosep} (3) already tells us that $\W(n)$ contains all genus 1 separating twists. 
It remains to show that $[T_a,f]\in\W(n)$ for  any $f\in \M_g^1$ and any separating curve $a$ of genus 2. When $g\ge5$, we have that 
$$[T_a,f]=T_aT_{f(a)}^{-1}\in \W(n)$$
by Proposition \ref{prop:difference of T2 in W}. 
If $g<5$, then $n=1$ or 2. When $n=1$, Johnson already proved that  $\W_g^1(1)=\Ch_g^1[2]=\I_g^1$ and hence there is nothing to be proved. When $n=2$ in which case $g$ must be 4, Proposition \ref{prop:BPtosep} implies that $\T(2)\le \W_4^1(2)$. 
\end{proof}

\begin{lemma}
\label{generating Ch[2n]}
    For $1\le n\le g-2$, the group $\Ch_g^1[2n]$ is normally generated by a genus $n$ bounding pair map $BP_n$, separating twists of genera 1 and 2, and an element $B_0:=T_aT_b^{-1}$ where $a,b$ are in Figure \ref{fig:B0} below.
    \begin{figure}[ht]
        \centering
        \includegraphics[width=0.6\linewidth]{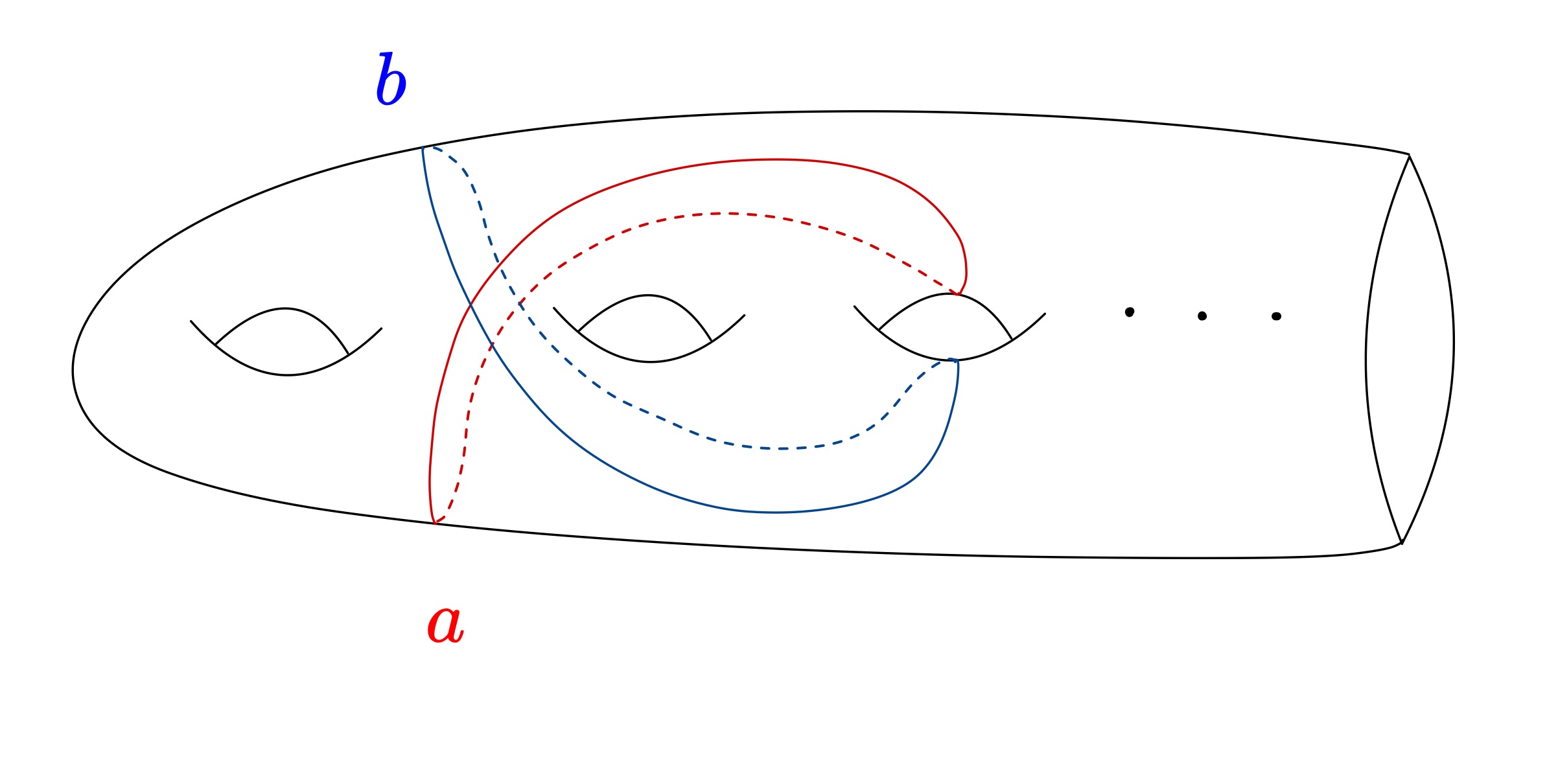}
        \vspace{-20pt}
        \caption{$B_0=T_aT_b^{-1}$ is called a homological genus 0 bounding pair map in \cite{Chilling}.}
        \label{fig:B0}
    \end{figure}

\end{lemma}
\begin{proof}
    There is a short exact sequence
    $$1\to \Ch_g^1\to \Ch_g^1[2n]\xrightarrow{\ch} 2nH\to1.$$
    Therefore, $\Ch_g^1[2n]$ is normally generated by $\Ch_g^1$ together with any element $f$ such that the $\Sp$-image of $e(f)$ generates $2nH$. By our computation in (\ref{ch of BPn}), we can take this $f$ to be $BP_n$. Next, $\Ch_g^1$ is normally generated by $B_0$ and $\J_g^1$ by Kosuge, \cite[Proposition 7] {Chilling}. $\J_g^1$ is generated by all separating twists of  genera 1 and 2 by Johnson, \cite[Theorem 1]{JohnsonLantern}. 
\end{proof}

\begin{theorem}
\label{main1}
When $1\le n\le g-2$, we have that
$$\W_g^1(n)=[\Ch_g^1[2n],\M_g^1] = [\W_g^1(n),\M_g^1].$$

\end{theorem}
\begin{proof}
    Since we always have that $\W(n)\le\Ch_g^1[2n]$ by (\ref{eq: W in Ch}), it suffices to prove that 
    $$[\Ch_g^1[2n],\M_g^1]\le \W(n)\le [\W(n),\M_g^1].$$
    We first show that $[\Ch_g^1[2n],\M_g^1]\le \W(n)$ by checking it on each of the normal generators of $\Ch_g^1[2n]$ from Lemma \ref{generating Ch[2n]}:
        $$\Ch_g^1[2n] =\langle\langle BP_n, B_0, T_1,T_2\rangle \rangle.$$
    By definition, $[BP_n,\M]\subseteq \W(n)$. Moreover, we can see that $B_0\in\W(n)$ by considering the curves $a,b,c$ in Figure \ref{fig:Bo in W} below.
    \begin{figure}[ht]
        \centering
        \includegraphics[width=0.6\linewidth]{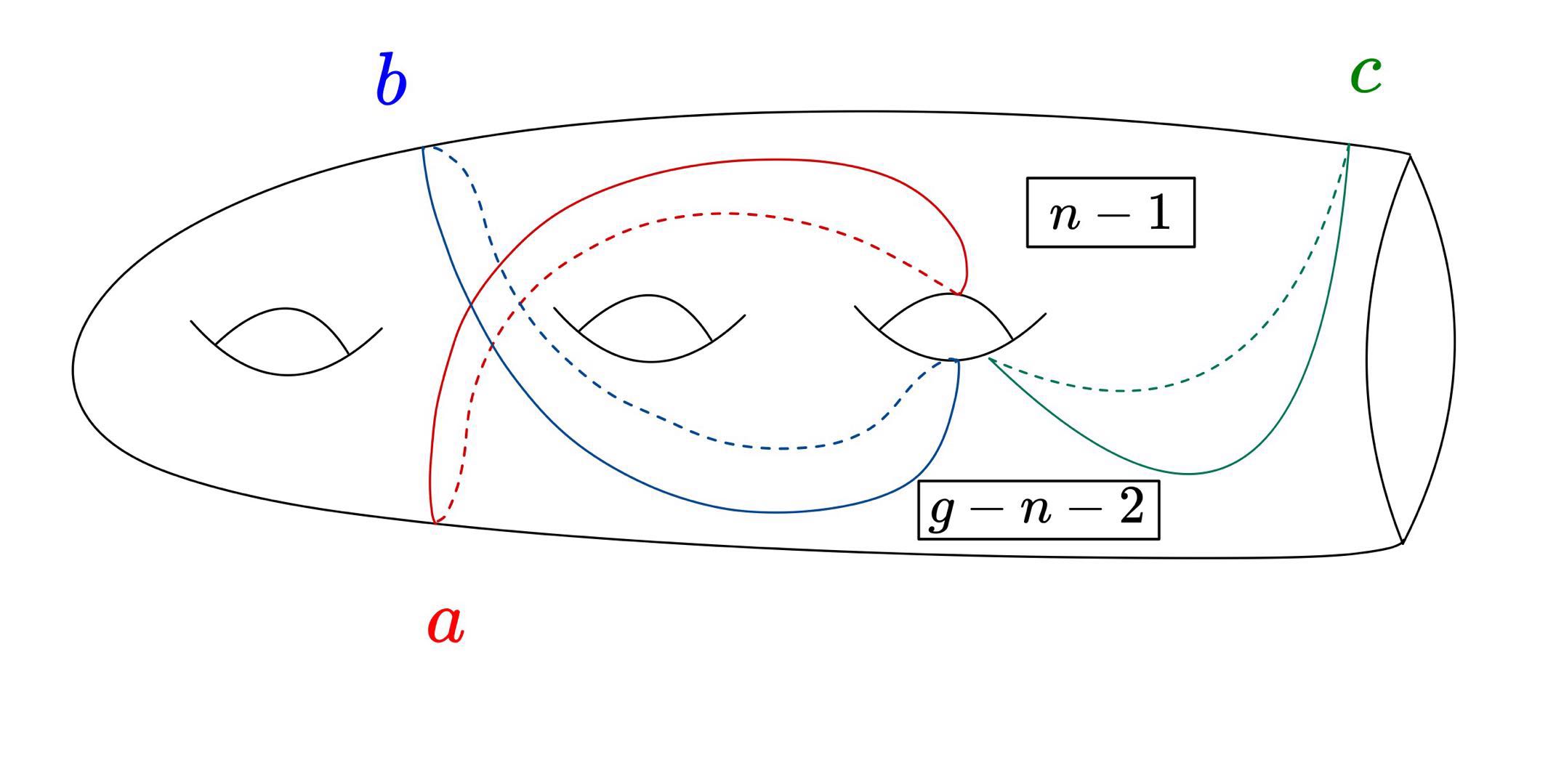}
                \vspace{-25pt}
        \caption{The numbers in the boxes represent the genera of subsurfaces that are not explicitly drawn. }
        \label{fig:Bo in W}
    \end{figure}
    Both $(a,c)$ and $(b,c)$ are bounding pairs of genus $n$. Therefore, we have that
    $$B_0=T_aT_b^{-1} = (T_aT_c^{-1})(T_bT_c^{-1})^{-1}\in\W(n).$$ 
    Since $1\le n\le g-2$, we have that $[\T(1),\M]\subseteq \T(1)\subseteq \W(n)$ by Proposition \ref{prop:BPtosep} and that $[\T(2),\M]\subseteq \W(n)$ by Proposition \ref{prop:difference of T2 in W}.

Next, we show that $\W(n)\le[\W(n),\M].$ Consider Figure \ref{fig:WM} where $(a,b)$ and $(c,d)$ are two genus $n$ bounding pairs such that $T_aT_b^{-1}(c)=d$. 

\begin{figure}[ht]
    \centering
    \includegraphics[width=0.7\linewidth]{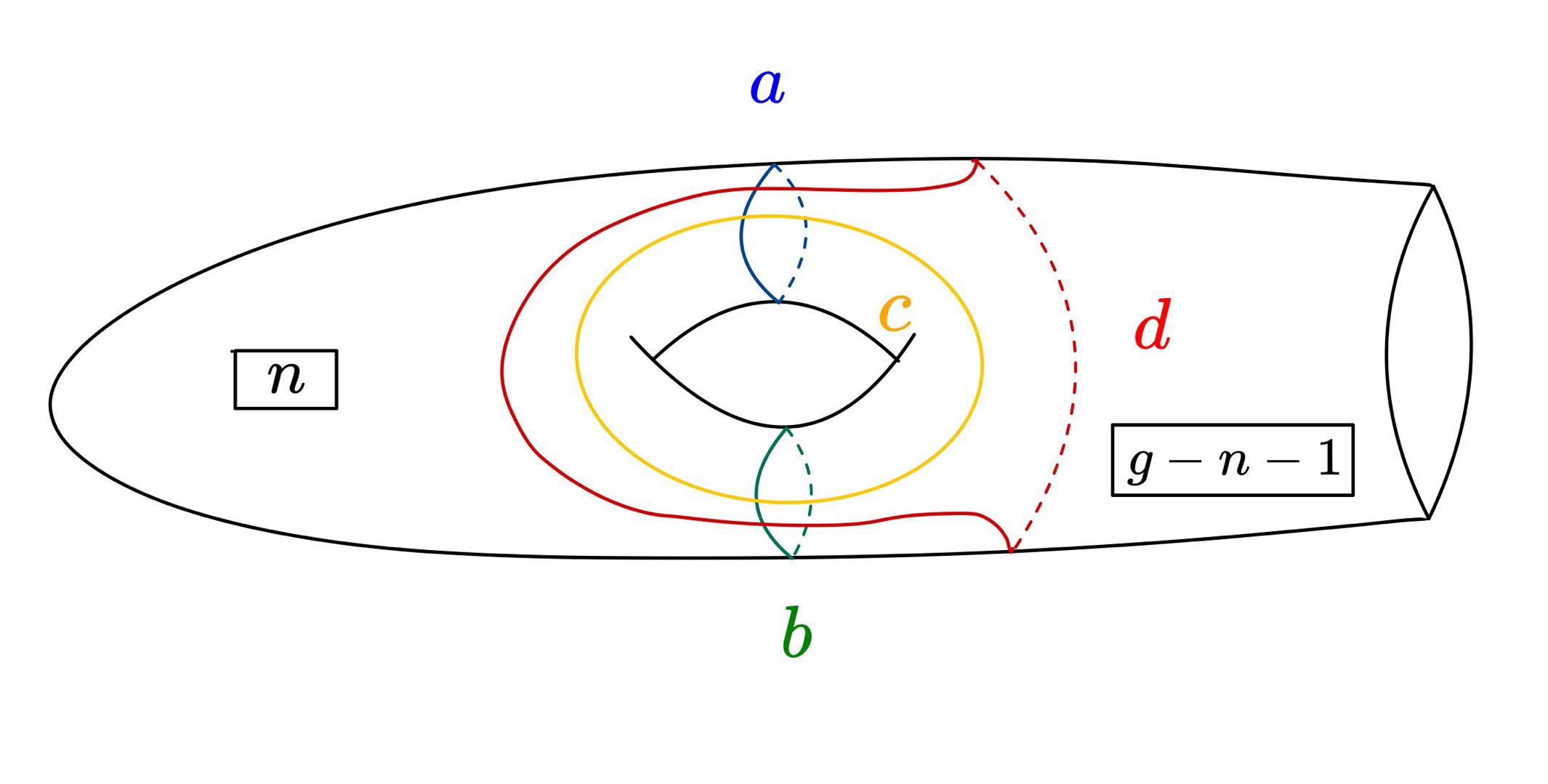}
    \vspace{-20 pt}
    \caption{$(a,b)$ and $(c,d)$ are two genus $n$ bounding pairs such that $T_aT_b^{-1}(c)=d$. The numbers in the boxes represent the genera of the subsurfaces.}
    \label{fig:WM}
\end{figure}
Hence, we have that 
$$T_cT_d^{-1}= [T_c, T_aT_b^{-1}]\in [\M,\W(n)].$$
This implies that $\W(n)\le[\W(n),\M].$

\end{proof}

\subsection{Proof of Theorem \ref{thm: W is kerd}  and Corollary \ref{cor:W vs K}}
Theorem \ref{thm: W is kerd} and Corollary \ref{cor:W vs K} will follow from Theorem \ref{main2} and Theorem \ref{index} below.
\begin{lemma}
\label{d on Tn Tn+1}
The image of $d$ on $\langle \T(n),\T(n+1)\rangle$ is $8n$ for $n$  odd and is $4n$ for $n$ even.
\end{lemma}
\begin{proof}
By Proposition \ref{d map properties}, we have that 
$d(T_n)=4n(n-1)$ and $d(T_{n+1})=4n(n+1)$. Hence,  the image of $d$ on $\langle \T(n),\T(n+1)\rangle$ is generated by $\gcd(4n(n-1),4n(n+1))$, which is $8n$ or $4n$ accordingly.
\end{proof}
\begin{theorem}
\label{main2}
If $1\le n\le g-2$, then  
$$\W_g^1(n)=\ker(d_{4n}:\Ch_g^1[2n]\to\Z/4n\Z).$$

\end{theorem}
\begin{proof}
Since $d_{4n}:\Ch_g^1[2n]\to \Z/4n\Z$ is a $\M_g^1$-invariant homomorphism by Lemma \ref{d4n invariant homo}, it maps the commutator subgroup $[\Ch_g^1[2n],\M_g^1]$ to zero. Hence, by Theorem \ref{main1}, we have that $\W(n)=[\Ch_g^1[2n],\M_g^1]\subseteq \ker d_{4n}.$ We just need to prove the reverse inclusion. 

Let $f\in \Ch_g^1[2n]$ be any element such that $d_{4n}(f)=0$ or equivalently $d(f)\in 4n\Z.$ By Proposition \ref{tau(W_n)=U2n}, we have that $\tau(\W(n))=U_{2n}=\tau(\Ch_g^1[2n])$. Therefore, there exists some $h\in\W(n)$ such that $\tau(f)=\tau(h)$. Since $fh^{-1}\in \ker\tau=\J_g^1$, we have that $d(fh^{-1})\in d(\J_g^1)=8\Z.$ In addition, by Proposition \ref{d map properties}, 
\[
d(fh^{-1})=d ( f ) + d ( h^{-1} ) +k( f ) \cdot k ( h^{-1} )
\]
is divisible by $4n$ since each of the three summands on the right-hand side is. Hence, $d(fh^{-1})$ is divisible by $\lcm(8,4n).$

In the case when $n$ is odd, $d(fh^{-1})$ is divisible by the least common multiple $\lcm(8,4n)=8n.$ By Lemma \ref{d on Tn Tn+1}, there exists an element $j\in \langle T_n,T_{n+1}\rangle$ such that $d(j)=d(fh^{-1})$. Then we have that 
\[d(fh^{-1}j^{-1})=0.\]
We conclude that $fh^{-1}j^{-1}\in \ker(d|_{\J_{g}^1})\le \W(n)$ by Lemma \ref{Kerd}. Since $h\in \W(n)$ and we know that $j\in \W(n)$ by Proposition \ref{prop:BPtosep}, it follows that $f\in \W(n).$

In the case when $n$ is even, $d(fh^{-1})$ is divisible by $\lcm(8,4n)=4n.$ By Lemma \ref{d on Tn Tn+1}, there still exists an element $j\in \langle T_n,T_{n+1}\rangle$ such that $d(j)=d(fh^{-1})$. The same proof applies. 
\end{proof}

Having determined the kernel of $d_{4n}$, we next compute its image. 
\begin{theorem}
\label{index}
For $1\le n\le g-2$, we have that 
$d_{4n}(\Ch_g^1[2n])=8\Z/4n\Z$ and that
\[
\frac{\Ch_g^1[2n]}{\W(n)}\cong\frac{\J_g^1}{\W(n)\cap \J_g^1} \cong\begin{cases}
			\Z/n\Z, & \text{if $n$ is odd}\\
            \Z/\frac{n}{2}\Z, & \text{if $n$ is even}.
		 \end{cases}
\]
Moreover, $\W(n)\cap \J_g^1$ is of infinite index in $\W(n)$.
\end{theorem}
\begin{proof}
The following commutative diagram relates the Johnson homomorphism $\tau$ in the vertical with the mod $4n$ Casson--Morita homomorphism $d_{4n}$ in the horizontal. For simplicity, we suppress the genus $g$ and boundary component from the notation. We have $\ker \tau=\J$ by Johnson's result and that $\ker d_{4n}=\W(n)$ by Theorem \ref{main1}.
\[\begin{tikzcd}
	& 0 & 0 \\
	0 & {\W(n)\cap\J} & \J & {d_{4n}(\J)} & 0 \\
	0 & {\W(n)} & {\Ch[2n]} & {d_{4n}(\Ch[2n])} & 0 \\
	& {\tau(\W(n))} & {\tau(\Ch[2n])} \\
	& 0 & 0
	\arrow[from=1-2, to=2-2]
	\arrow[from=1-3, to=2-3]
	\arrow[from=2-1, to=2-2]
	\arrow[from=2-2, to=2-3]
	\arrow[from=2-2, to=3-2]
	\arrow["{d_{4n}}", from=2-3, to=2-4]
	\arrow[from=2-3, to=3-3]
	\arrow[from=2-4, to=2-5]
	\arrow["j", from=2-4, to=3-4]
	\arrow[from=3-1, to=3-2]
	\arrow[from=3-2, to=3-3]
	\arrow["\tau", from=3-2, to=4-2]
	\arrow["{d_{4n}}", from=3-3, to=3-4]
	\arrow["\tau", from=3-3, to=4-3]
	\arrow[from=3-4, to=3-5]
	\arrow["i", "="', from=4-2, to=4-3]
	\arrow[from=4-2, to=5-2]
	\arrow[from=4-3, to=5-3]
\end{tikzcd}\]
By Proposition \ref{tau(W_n)=U2n} we have that $\tau(\W(n))=U_{2n} = \tau (\Ch[2n])$. As a consequence, $\W(n)\cap \J$ is of infinite index in $\W(n)$. Moreover, the bottom-left horizontal map $i$ is an equality. By the Snake Lemma or by diagram-chasing, the top-right vertical map $j$ is also an equality. Hence, we have 
$$\frac{\Ch[2n]}{\W(n)}\cong d_{4n}(\Ch[2n])\overset{j}{=} d_{4n}(\J)\cong \frac{\J}{\W(n)\cap \J}.$$
Next, we will determine the isomorphic groups above by computing $d_{4n}(\J)$. Proposition \ref{d map properties} and Johnson's generation result for $\J$ imply that $d(\J)=8\Z.$ Therefore, we have that 
$$d_{4n}(\J) =8(\Z/4n\Z)\cong\begin{cases}
			\Z/n\Z, & \text{if $n$ is odd}\\
            \Z/\frac{n}{2}\Z, & \text{if $n$ is even}
		 \end{cases}$$
where $8(\Z/4n\Z)$ represents the subgroup of $\Z/4n\Z$ generated by $8$ mod $4n$. 

\end{proof}

\subsection{Proof of Corollary \ref{cor:homological results}}
To prove (1), we have that 
\begin{align*}
    H_1(\Ch_g^1[2n];\Z)_{\M_g^1} &\cong \frac{\Ch_g^1[2n]}{[\Ch_g^1[2n],\M_g^1]}\cong \frac{\Ch_g^1[2n]}{\ker(d_{4n})} &\text{by Theorem \ref{main1} and \ref{main2}}\\
    &\cong \im(d_{4n})\cong \begin{cases}
			\Z/n\Z, & \text{if $n$ is odd}\\
            \Z/\frac{n}{2}\Z, & \text{if $n$ is even}
		 \end{cases}
         &\text{by Theorem \ref{index}.}
\end{align*}
To prove (2), we have that by Theorem \ref{main1}
$$H_1(\W_g^1(n);\Z)_{\M_g^1} \cong \frac{\W_g^1(n)}{[\W_g^1(n),\M_g^1]}=0.$$

\subsection{Image of $d$ on various subgroups of Torelli groups}
\label{sec: integral d}
In the proof of Theorem \ref{index} above, we computed the image of $d_{4n}=(d\mod 4n)$ on $\Ch_g^1[2n]$ and $\W_g^1(n)$. In this subsection, we will use those results to calculate the image of the integral $d$ on various subgroups of Torelli groups.

\begin{proposition}
    \label{d(BPn)}
    Consider the genus $n$ bounding pair map $T_cT_d^{-1}$ pictured in Figure \ref{fig:WM}. We have that 
    $$d(T_cT_d^{-1})=4n(1-4n).$$
\end{proposition}
\begin{proof}
    For any  $h\in \I_g^1$ and $f\in \M_{g,1}$, we have that 
\[
k(fhf^{-1})=f_*k(fh)+k(f^{-1})=f_*k(f)+f_*k(h)-f_*k(f)=f_*k(h)
\]
\begin{equation}
\begin{aligned}
d(fhf^{-1}h^{-1}) &= d(fhf^{-1})+d(h^{-1})+k(fhf^{-1})\cdot k(h^{-1}) \\
&= d(h)+k(f)\cdot (k(h)+k(fh))+d(h^{-1})+f_*k(h)\cdot k(h^{-1}) \\
&= k(f)\cdot (k(h)+k(f)+k(h))-f_*k(h)\cdot k(h)\\
&=2k(f)\cdot k(h)-f_*k(h)\cdot k(h).
\end{aligned}
\end{equation}
In particular, let $a,b,c,d$ denote the curves in Figure \ref{fig:WM} and let $h:=T_aT_b^{-1}$ be a genus $n$ bounding pair map. We have that 
$$d(T_cT_d^{-1})=d(T_chT_c^{-1}h^{-1}) = 2k(T_c)\cdot k(h)-(T_c)_*k(h)\cdot k(h).$$
By (\ref{ch of BPn}), we have that $k(h)=2nb_{n+1} = 2na.$

Next, we calculate $k(T_c)$ using  Morita's definition of $k$ which we described in the second paragraph of Section \ref{sec:Def of d}.  We fix based loops $\alpha_i, \beta_i$ for $i=1,\dots,g$ as generators for $\pi_1(S_g^1)$ as in Figure \ref{fig:generators} below.
\begin{figure}[ht]
    \centering
    \includegraphics[width=0.65\linewidth]{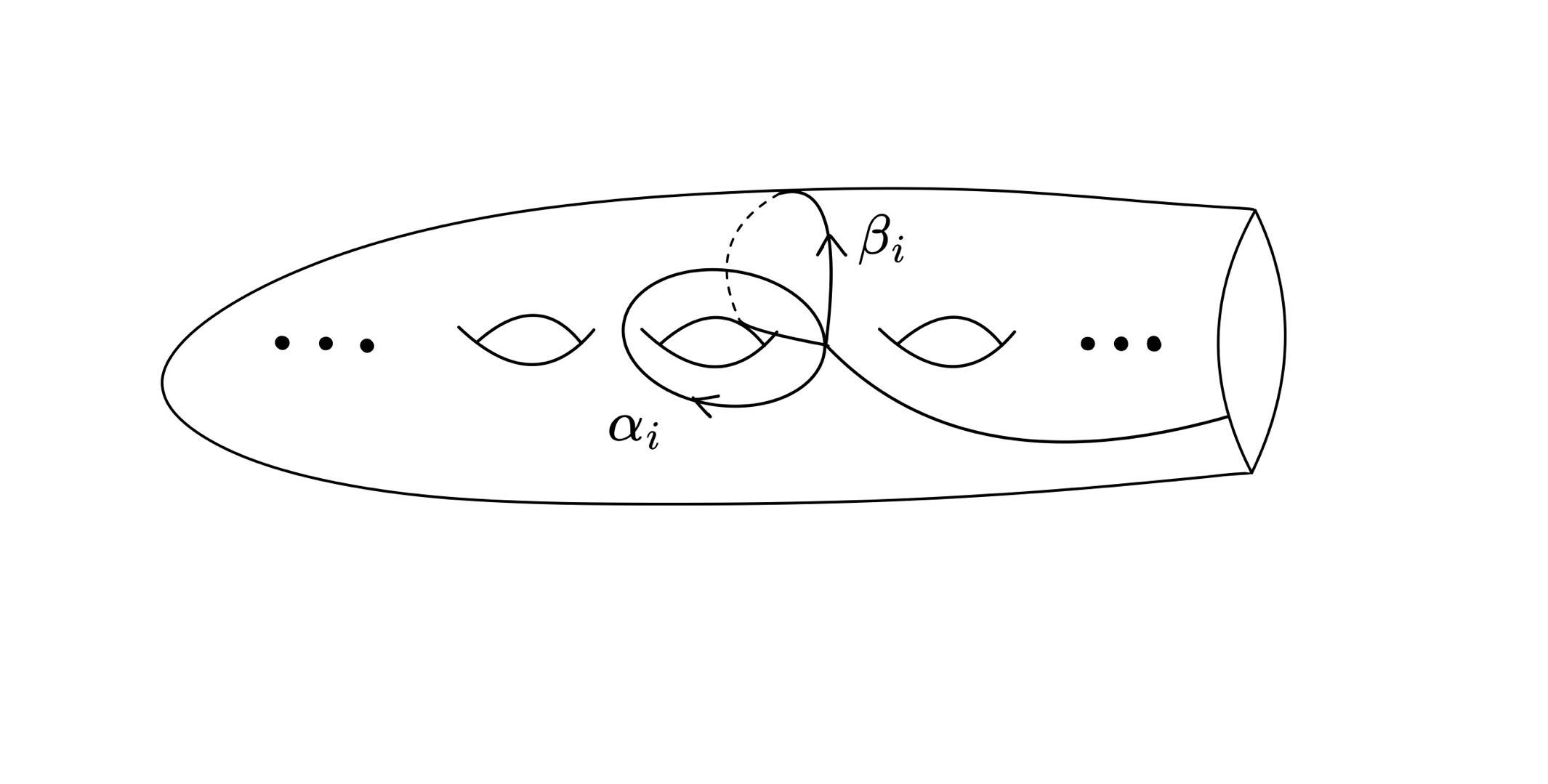}
    \vspace{-35 pt}
    \caption{$\alpha_i, \beta_i$ as generators for $\pi_1(S_g^1)$}
    \label{fig:generators}
\end{figure}

We have that $T_c$ fixes all the generators except that $T_c(\beta_{n+1})=\beta_{n+1}\alpha_{n+1}^{-1}.$ Hence, $k^\vee(T_c)$ takes all the generators to zero except that 
$$k^\vee(T_c)(\beta_{n+1})= \epsilon(\beta_{n+1}\alpha_{n+1}^{-1})-\epsilon(\beta_{n+1}) = [\beta_{n+1}]\cdot[\alpha_{n+1}^{-1}] = 1.$$
This implies that $k(T_c)=c.$ 
Finally, we conclude that 
$$d(T_cT_d^{-1})=2k(T_c)\cdot k(h)-(T_c)_*k(h)\cdot k(h)=(2c)\cdot (2na)-(T_c)_*(2na)\cdot (2na)=4n-4n^2.$$
\end{proof}

\begin{proposition}
    \label{prop:dCh[2n]}
    For $1\le n\le g-2$, we have that  
    \begin{enumerate}
        \item $d(\Ch_g^1[2n])=\begin{cases}
            8\Z& \text{if $n$ is even}\\
            4\Z& \text{if $n$ is odd},
        \end{cases}$
        \item $d(\W_g^1(n))=4n\Z.$
    \end{enumerate} 
\end{proposition}
One immediate corollary is that both $d(\Ch_g^1[2n])$ and $d(\W_g^1(n))$ are  subgroups of $\Z$. This is not obvious because $d$ is not a group homomorphism on either $\Ch_g^1[2n]$ or $\W_g^1(n)$. 

\begin{proof}
We first consider (1). Theorem \ref{index} implies that 
    $$d(\Ch_g^1[2n])\Big({\Z}/{4n\Z}\Big) = d_{4n}(\Ch_g^1[2n]) 
    \overset{j}{\cong} d_{4n}(\J_g^1)=8({\Z}/{4n\Z}).
$$
where the first group denotes the subgroup of $\Z/4n\Z$ generated by the subset $d(\Ch_g^1[2n])$ and where $j$ is the same isomorphism as in the proof of Theorem \ref{index}. The equalities above imply that 
$$8\Z\subseteq d(\Ch_g^1[2n])\subseteq\gcd(8,4n)\Z$$
If $n$ is even, then the inclusions above imply that $d(\Ch_g^1[2n])=8\Z$. Suppose $n$ is odd. We have that $8\Z\subseteq d(\Ch_g^1[2n])\subseteq 4\Z$. Hence, in order to prove that $d(\Ch_g^1[2n])= 4\Z$, it remains to show that every integer of the form $4+8m$ is in $d(\Ch_g^1[2n])$. Let $h$ denote the genus $n$ bounding pair map $h\in\Ch_g^1[2n]$ as in Proposition \ref{d(BPn)} such that $d(h)=4n(1-4n)$. Since $n$ is odd, $d(h)=4 \pmod 8$ and hence there exists $r\in\Z$ such that 
$4+8m=d(h)+8r.$
Let $T_2$ denote any genus 2 separating twist. We take an element $hT_2^r\in \Ch_g^1[2n]$ which satisfies that
\begin{align*}
    d(hT_2^r)&= d(h)+d(T_2^r)+k(h)\cdot (T_2^r)_*k(T_2^r)+0 &\text{by Proposition \ref{d map properties}}\\
    &=d(h)+d(T_2^r)&\text{since $k(\J_g^1)=0$}\\
    &=d(h)+rd(T_2)&\text{since $d$ is a homomorphism on $\J_g^1$}\\
    &= d(h)+8r=4+8m.
\end{align*}
Hence, we conclude that $d(\Ch_g^1[2n])= 4\Z$. 

We will prove (2) using a similar argument. Theorem \ref{main2} implies that 
$$d(\W(n))\subseteq 4n\Z$$
If $n$ is even, then by Proposition \ref{prop:BPtosep} and Proposition \ref{d on Tn Tn+1}, we have  that 
$$4n\Z= d\big(\langle \T(n),\T(n+1)\rangle\big)\subseteq d( \W(n)).$$
We are done in this case. 

If $n$ is odd, then by Proposition \ref{prop:BPtosep} and Proposition \ref{d on Tn Tn+1}, we have  that 
$$8n\Z= d\big(\langle \T(n),\T(n+1)\rangle\big)\subseteq d( \W(n)).$$
In order to prove that $d( \W(n))=4n\Z$, it remains to show that every integer of the form $4n+8nm$ is in $d( \W(n))$. Let $h$ again denote the genus $n$ bounding pair map as above such that $d(h)=4n(1-4n)$. There exists an integer $r$ such that 
$$4n+8nm=4n(1-4n)+8r.$$
We take an element $f\in \langle \T(n),\T(n+1)\rangle$ such that $d(f)=8n$ and an element $hf^r\in \W_g^1(n)$. We have that
\begin{align*}
    d(hf^r)&= d(h)+d(f^r)+k(h)\cdot f_*k(f^r)+0 &\text{by Proposition \ref{d map properties}}\\
    &=d(h)+d(f^r)&\text{since $k(\J_g^1)=0$}\\
    &=d(h)+rd(f) = d(h)+8r=4n+8nm.
\end{align*}
Hence, we conclude that $d(\W_g^1(n))= 4n\Z$. 
\end{proof}

\begin{corollary}
    For $g\ge 3$, we have that $d(\I_g^1)= 4\Z$. 
\end{corollary}
\begin{proof}
    $\I_g^1=\W_g^1(1)=\Ch_g^1[2]$. 
\end{proof}

\section{Punctured and closed surfaces}
\label{sec:closed and punctured}

\subsection{Point-pushing subgroups and genus $g-1$ bounding pair maps}
Let $S_g^1$, $S_{g,1}$ and $S_g$ denote a once-bordered surface, a once-punctured surface, and a closed surface, all of genus $g$, respectively.  Their Torelli groups are related by the following commutative diagram of exact sequences when $g\ge 2$:
\begin{equation}
    \label{big commutative diagram}
    \begin{tikzcd}
 & 0 \arrow{d} &   0 \arrow{d} & &\\
 & \mathbb{Z} \arrow{d} \arrow{r}{\cong} &  \mathbb{Z}=\langle T_b\rangle \arrow{d} & &\\
0 \arrow{r} & \pi_{1}(UT(S_g))\arrow{d} \arrow{r}{Push} & \I_g^1 \arrow{r}{Cap} \arrow{d}{cap} & \I_g \arrow{r} \arrow{d}{\cong} & 0 \\
0 \arrow{r} & \pi_{1}(S_g) \arrow{r}{push} \arrow{d} & \I_{g,1} \arrow{r}{forget} \arrow{d} & \I_g \arrow{r} & 0 \\
 & 0 & 0 & &
\end{tikzcd}
\end{equation}
Here $T_b$ denotes the Dehn twist along a curve isotopic to the boundary. The map $Cap$ (or $cap$) corresponds to attaching a disk (or a once-punctured disk) to the boundary and extending the mapping class by the identity to the disk.  The last row is the Birman exact sequence where $push$ denotes the point-pushing map and $forget$ denotes forgetting the marked point. See \emph{e.g.} Section 4.2.5 in \cite{primer} for additional details.

\begin{proposition}\label{g-1 case}
When $g\ge 2$, we have that
\begin{enumerate}
    \item $\W_g^1(g-1)=Push(\pi_1(UT(S_g)))$.
    \item $\W_{g,1}(g-1)=push(\pi_1(S_g))$.
\end{enumerate}
\end{proposition}
\noindent Note that for the closed surfaces, we have that $\W_g(g-1)=0$.
\begin{proof}
     First we consider part (2). For a based simple closed curve $a$ in $S_g$, we know that $push(a) = T_cT_d^{-1}$ where $c,d$ are boundary components of a tubular neighborhood of $a$. In particular, $push(a)$ is a genus $g-1$ bounding pair map. Conversely, for any genus $g-1$ bounding pair map $T_cT_d^{-1}$, the bounding pair $(c,d)$ cuts $S_g$ into two subsurfaces, one of which is a pair of pants $S_0^3$. Hence, we can find a based simple closed curve $a$ such that $push(a)=T_cT_d^{-1}$. Part (2) is proved.

     To prove (1), it suffices to prove that $T_b\in \W_g^1(g-1)$, which was established in Proposition \ref{prop:BPtosep}, part (2). 
\end{proof}

\subsection{$\W(n)$ for once-punctured surfaces}
Consider the second vertical column in (\ref{big commutative diagram}):
\begin{equation}
    \label{cap SES}
    1\to \langle T_b\rangle \to \I_g^1\to \I_{g,1}\to 1.
\end{equation}
Since the Johnson homomorphism $\tau_g^1$ vanishes on $T_b$, it descends to a map $\tau_{g,1}$ on the quotient $\I_{g,1}$, and so does the Chillingworth homomorphism $\ch=C\circ \tau$. For any $n\in\Z_{>0}$, we define
$$\Ch_{g,1}[2n]:= \{f\in\I_{g,1}\ :\ C(\tau_{g,1}(f))=0 \pmod {2n}\}.$$
Let $\W_{g,1}(n)$ be the group generated by all genus $n$ bounding pair maps of  $S_{g,1}$. 
\begin{theorem}
\label{punctured  W}
For $1\le n\le g-2$ and for $\ell = \gcd(g(g-1),n)$, 
there is a well-defined homomorphism
$d_{4\ell}:= (d\mod 4\ell)$
such that
\[
\W_{g,1}(n)=\ker(d_{4\ell}: \Ch_{g,1}[2n]\to\Z/4\ell\Z).
\]
Moreover, in this case, we have that  $d_{4\ell}(\Ch_{g,1}[2n])=8\Z/4\ell\Z$ and hence
    \[
\frac{\Ch_{g,1}[2n]}{\W_{g,1}(n)}\cong\frac{\J_{g,1}}{\W_{g,1}(n)\cap \J_{g,1}} \cong\begin{cases}
			\Z/\ell\Z, & \text{if $n$ is odd}\\
            \Z/\frac{\ell}{2}\Z, & \text{if $n$ is even}.
		 \end{cases}
\]
\end{theorem}

\begin{proof}
The exact sequence   (\ref{cap SES}) restricted to the subgroup $\Ch_g^1[2n]$ gives
\[
1\to\langle T_b\rangle \to \Ch_g^1[2n]\to \Ch_{g,1}[2n]\to 1.
\]
Moreover, the surjective capping homomorphism takes $\W_g^1(n)$ onto $\W_{g,1}(n)$.

Theorem \ref{main2} tells us that $\W_g^1(n)$ is the kernel of 
\[\Ch_g^1[2n]\xrightarrow{d_{4n}} \Z/4n\Z.\]
Since $b$ is a separating curve of genus $g$, we have that $d(T_b)=4g(g-1)$ by Proposition \ref{d map properties}. Since $\ell = \gcd(g(g-1),n)$, the map $d_{4\ell}$ descends to a well-defined map on the quotient $\Ch_{g,1}[2n]$ whose kernel is exactly $\W_{g,1}(n)$. The calculation for the image of $d_{4\ell}$ is similar as in the proof of Theorem \ref{index}.
\end{proof}

\subsection{$\W(n)$ for closed surfaces}
Consider the Birman exact sequence: 
\begin{equation}
    \label{birman SES}
1\to \pi_1(S_g)\xrightarrow{push} \I_{g,1}\to \I_g\to 1
\end{equation}
The Johnson homomorphism $\tau_{g,1}$ descends to a well-defined map $\tau_g$ on the quotient
\[
\tau_g: \I_g\to \wedge^3 H/H
\]
where the inclusion of $H$ into $\wedge^3 H$ is given by $x\mapsto \sum_{i=1}^g a_i\wedge b_i\wedge x$ where $a_i,b_i$ represent a symplectic basis for $H$. Similarly, the Chillingworth homomorphism is the composition
$$\I_g\xrightarrow{\tau_g} \wedge^3 H/H \xrightarrow{C} H/(2g-2)H.$$
For any $n\in\Z_{>0}$ we define
$$\Ch_{g}[2n]:= \{f\in\I_{g}\ :\ C(\tau_{g}(f))=0 \pmod{2n}\}.$$
From the definition we can see that $\Ch_g[2n]=\Ch_g[2m]$ if $m=\gcd(g-1,n)$. 
Let $\W_{g}(n)$ be the group generated by all genus $n$ bounding pair maps of $S_{g}$.

\begin{theorem}
\label{thm: closed W}
For $1\le n\le g-2$ and for   $m:=\gcd(g-1,n)$, we have that
\begin{equation}
    \label{closed Wn for gcd}
    \W_g(n)=\W_g(m).
\end{equation}
Moreover, there is a well-defined homomorphism
$d_{4m}:=(d \mod 4m)$ such that 
\[
\W_{g}(n)=\ker(d_{4m}: \Ch_{g}[2n]\to\Z/4m\Z).
\]
In this case, we have that  $d_{4m}(\Ch_{g}[2n])=8\Z/4m\Z$ and 
\[
\frac{\Ch_{g}[2n]}{\W_{g}(n)}\cong\frac{\J_{g}}{\W_{g}(n)\cap \J_{g}} \cong\begin{cases}
			\Z/m\Z, & \text{if $m$ is odd}\\
            \Z/\frac{m}{2}\Z, & \text{if $m$ is even}.
		 \end{cases}
\]

\end{theorem}
\begin{proof}
    First of all, to prove (\ref{closed Wn for gcd}), we first make the following observations:
    \begin{enumerate}
\item
For any integer $m$, if $\W_g(m)\subseteq \W_g(n)$ and $m+n\le g-1$, then $\W(n+m)\subseteq \W(n)$. To see this, consider Figure \ref{fig:abc} below. We have that $BP_{n+m}=T_aT_c^{-1} = (T_aT_b^{-1})(T_bT_c^{-1})$ which is in $\W_g(n)$.
\begin{figure}[b]
    \centering
    \includegraphics[width=0.6\linewidth]{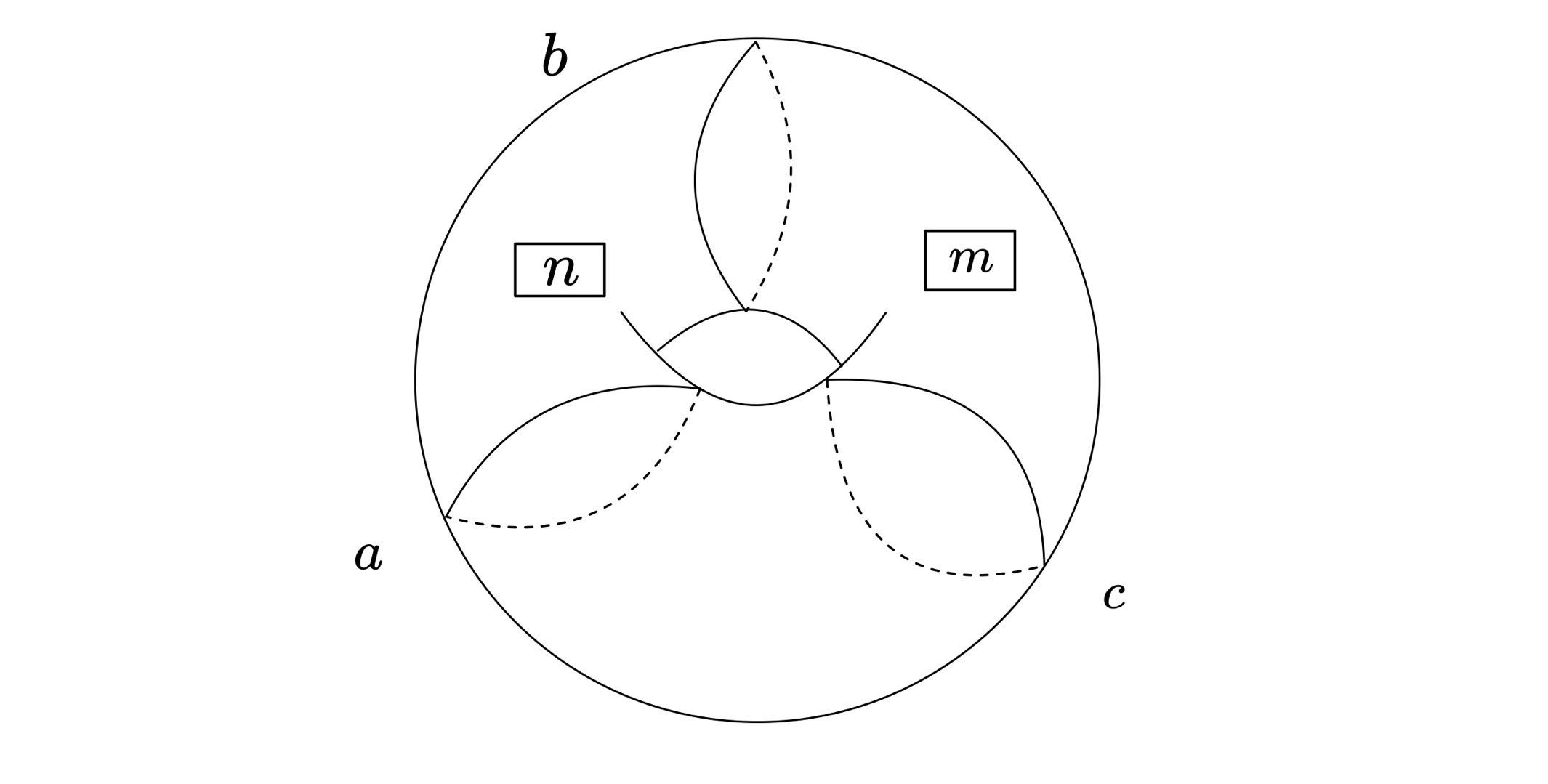}
    \vspace{-10pt}
    \caption{The numbers in the boxes denote the genera of the subsurfaces.}
    \label{fig:abc}
\end{figure}
\item 
For any integer $m$, if $W(m)\subseteq \W(n)$ and $m+n> g-1$, we have that $\W(n+m-g+1)\subseteq \W(n)$. The proof is similar to the proof of (1) above. 
\end{enumerate}
Observation (1) implies that $\W(n)\subseteq \W({\gcd(n,g-1)})$. The reverse inclusion follows from both (1) and (2) together with the Euclidean algorithm. 

Next, we show that $\W_g(n)=\ker d_{4m}$ where $m=\gcd(g-1,n)$. By (\ref{closed Wn for gcd}), without loss of generality we can assume that $n=m$ or equivalently that $n$ divides $g-1$. We first need to check that $d_{4n}$ gives a well-defined map on $\Ch_g[2n].$ Consider the Birman exact sequence (\ref{birman SES}) where the surjective forget-the-marked-point homomorphism takes the subgroup $\W_{g,1}(n)$ onto $\W_{g}(n)$.
If $\alpha$ is a nonseparating simple closed curve, then $push(\alpha)$ is a bounding pair map of genus $g-1$. Hence, we conclude that its image under the Chillingworth homomorphism  (\ref{ch of BPn}) is
$$\ch(push(\alpha))\in 2(g-1)H.$$
Since $g-1$ is divisible by $n$, we have that $\ch(push(\alpha))\in 2nH$. Therefore, 
$$push(\pi_1(S_g))\subseteq \ch^{-1}(2nH)=\Ch_{g,1}[2n].$$ This gives an exact sequence
\[
1\to \pi_1(S_g)\xrightarrow{push} \Ch_{g,1}[2n]\to \Ch_g[2n]\to 1.
\]
Since $g-1$ is divisible by $n$, the homomorphism in Theorem \ref{punctured  W} is $d_{4l}=d_{4n}$
\begin{equation}
    \label{d4n for punctured}
    d_{4n}:\Ch_{g,1}[2n]\to \Z/4n\Z.
\end{equation}
Moreover, Theorem \ref{punctured  W} gives that
\begin{align*}
    \W_{g,1}(g-1) &= \ker (\Ch_{g,1}[2(g-1)]\xrightarrow{d_{4(g-1)}} \Z/4(g-1)\Z)\\
    &\subseteq \ker (\Ch_{g,1}[2n]\xrightarrow{d_{4n}} \Z/4n\Z) = \W_{g,1}(n)
\end{align*}
In particular, for any simple closed curve $\alpha$, we have that $push(\alpha)\in \W_{g,1}(n)$. Therefore, the homomorphism $d_{4n}$ in (\ref{d4n for punctured}) descends to a well-defined map on the quotient
$$d_{4n}:\Ch_g[2n]\to \Z/4n\Z$$
whose kernel is exactly $\W_g(n).$ The calculation for image of $d_{4n}$ is exactly the same as in the once-bordered and once-punctured cases already proven.
\end{proof}

\begin{corollary}
\label{prop:warmup}
For $1 \leq n \le g-2$, we have that  $\W_g(n)= \I_g$ if and only if $n$ is relatively prime to $g-1$.
\end{corollary}
\begin{proof}
    Let $m:=\gcd(g-1,n)$. By (\ref{closed Wn for gcd}), we have that     $\W_g(n)=\W_g(m).$ If $n$ is relatively prime to $g-1$, then $m=1$, and $\W_g(m)=\W_g(1)=\I_g$ since $g\ge 3.$

    Conversely, suppose that $m\ge2$. Consider a bounding pair map of genus $m-1\ge1$ which is an element in $\I_g$. By (\ref{ch of BPn}), we have that
    $$e(BP_{m-1})=2(m-1)b_m \ \ \ \text{ and } \ \ \ e(BP_{m})=2mb_{m+1}.$$
    Since  $e(BP_{m-1})\notin e(\W_g(m))$, we have that $BP_{m-1}\notin \W_g(m)$ and $\W_g(n) \neq \I_g$. 
\end{proof}

\bibliographystyle{plain}
\bibliography{BPsepV2}

\end{document}